\documentclass[amsmath,amssymb,showkeys, showpacs, superscriptaddress]{revtex4}
\usepackage{graphicx}
\usepackage{graphicx}
\usepackage{amsthm}
\usepackage{amsmath}
\usepackage{amssymb}
\usepackage{color}

\expandafter\ifx\csname package@font\endcsname\relax\else
 \expandafter\expandafter
 \expandafter\usepackage
 \expandafter\expandafter
 \expandafter{\csname package@font\endcsname}%
\fi

\theoremstyle{definition}

\parskip \baselineskip

\begin{document}
\vspace*{0.4 in}

\title{On W. Gordon's integral (1929) and related identities}

\author{Nasser Saad}
\email{nsaad@upei.ca}  
\affiliation{Department of Mathematics and Statistics, University of Prince Edward Island,\\ Charlottetown, Prince Edward Island  C1A 4P3, Canada.
}
\def\dbox#1{\hbox{\vrule  
                        \vbox{\hrule \vskip #1
                             \hbox{\hskip #1
                                 \vbox{\hsize=#1}%
                              \hskip #1}%
                         \vskip #1 \hrule}%
                      \vrule}}
\def\qed{\hfill \dbox{0.05true in}}  
\def\square{\dbox{0.02true in}} 
\begin{abstract}
\noindent Analytic evaluation of Gordon's integral 
$$\operatorname {J}_c^{j(\pm p)}(b,b';\lambda,w,z)=\int_0^\infty x^{c+j-1}e^{-\lambda x}{}_1F_1(b;c;wx){}_1F_1(b';c\pm p;zx)dx,$$
are given along with convergence conditions. It shows enormous number of definite integrals, frequently appear in theoretical and mathematical physics applications, easily deduced from this generalized integral.
\end{abstract}

\maketitle
\section{Gordon's integral: Introduction}
\noindent Among the important integrals in theoretical and mathematical physics is W. Gordon's integral \cite{gordon}, see also \cite{karule,landau,saad2003,tarasov2003},
\begin{align}\label{integral1}
\operatorname {J}_c^{j(\pm p)}(b,b';\lambda,w,z)&=\int_0^\infty x^{c+j-1}e^{-\lambda x}{}_1F_1(b;c;wx){}_1F_1(b';c\pm p;zx)dx\notag\\
&(c+j>0;~\lambda>0;~c,c\pm p\neq 0,-1,-2,\dots;~p\geq 0;~j=0,\pm 1,\pm 2,\dots),
\end{align}
where ${}_1F_1$ is the confluent hypergeometric function ${}_1F_1(b;c;z)=\sum_{k=0}^\infty{(b)_k\,z^k}/[(c)_k\,{k!}]$ in which $(b)_k=b(b+1)\dots(b+n-1)=\Gamma(b+k)/\Gamma(b)$ is the Pochhammer symbol defined in terms of Gamma function. The massive uses of this integral and the subclasses of it span large volume of research papers and monographs  \cite{karule,landau,opps2009,saad2003,tarasov2003}.
It was proven (Lemma 1 in \cite{saad2003}) that, for  $c+j>0$ and $|w|+|z|<|\lambda|$, 
\begin{align}\label{integral2}
\operatorname {J}_c^{j(\pm p)}(b,b';\lambda,w,z)&
=\dfrac{\Gamma(c+j)}{\lambda^{c+j}}\,F_2\left(\begin{matrix}c+j;&b,&b\rq{}\\ 
~&c,&c\pm p\end{matrix};\frac{w}{\lambda},\frac{z}{\lambda}\right),
\end{align} 
where the second Appell function reads (\cite{appell}, equation (2))  
\begin{align}\label{def1}
\setlength\arraycolsep{1pt}
F_2\left(\begin{matrix}a;&b,&b'\\
~&c,&c'\end{matrix};w,z\right)&\equiv F_{2}\left(  a; b,b^{\prime};c,c^{\prime};w,z\right)  =\sum_{m=0}^{\infty}\sum
_{p=0}^{\infty}\frac{\left(  a\right)_{m+p}\left(b\right)_{m}\left(
b^{\prime}\right)_{p}}{\left(c\right)_{m}\left(c^\prime\right)_{p}}\frac{w^{m}\,z^{p}}{m!\,p!},\quad (c,c^\prime \neq 0,-1,\dots; |w|+ |z| <1).
\end{align}
Exact analytical expressions of this integral by means of more elementary functions are given in the present work, where many subclasses are analysed and evaluated in simplified expressions allow for faster computations.

\section{Gordon's integral: Closed form expressions}
\noindent By means of the double integral representation of the second Appell function (\cite{appell}, equation 7; see also \cite{saad2003}), for $c,c-p\neq 0,-1,-2,\dots$, $j,p=0,1,2,\dots,~|w|+|z|<1$, it follows, for $j\geq p$, that
\begin{align}\label{integral4}
\operatorname {J}_c^{j(\pm p)}(b,b';\lambda,w,z)&=\int_0^\infty x^{c+j-1}e^{-\lambda x}{}_1F_1(b;c;wx){}_1F_1(b';c\pm p;zx)dx\notag\\
&=\dfrac{\Gamma(c+j)}{\lambda^{c+j-b'}(\lambda-z)^{b'}}\sum_{k=0}^{j\mp p} \dfrac{(-j\pm p)_k(b')_k}{(c\pm p)_k~k!}\left(1-\dfrac{\lambda}{z}\right)^{-k}F_1\left(b,c+j-b',b'+k;c;\dfrac{w}{\lambda},\dfrac{w}{\lambda-z}\right),\notag\\
&(c+j>0;~\lambda>0;~c,c\pm p\neq 0,-1,\dots;~p\geq 0;~j=0,\pm 1,\dots;~|w|+|z|<\lambda),
\end{align}
particularly, for $p=j= 0,1,2,\dots$,
\begin{align}\label{integral5}
\operatorname {J}_c^{jj}(b,b';\lambda,w,z)&=\int_0^\infty x^{c+j-1}e^{-\lambda x}{}_1F_1(b;c;wx){}_1F_1(b';c+ j;zx)dx=\dfrac{\Gamma(c+j)\,F_1\left(b,c+ j-b',b';c;\dfrac{w}{\lambda},\dfrac{w}{\lambda-z}\right)}{\lambda^{c+j-b'}(\lambda-z)^{b'}},\notag\\
&(c+j>0;~\lambda>0;~c\neq 0,-1,\dots;~j=0,\pm 1,\dots;~|w|+|z|<\lambda),\notag\\
\operatorname {J}_c^{jj}(c+j,b';\lambda,w,z)&=\int_0^\infty x^{c+j-1}e^{-\lambda x}{}_1F_1(c+j;c;wx){}_1F_1(b';c+ j;zx)dx\notag\\
&=\dfrac{\Gamma(c+j)(\lambda-w)^{b'-c-j}}{(\lambda-z-w)^{b'}}F_1\left(-j,c+ j-b',b';c;\dfrac{w}{w-\lambda},\dfrac{w}{w+z-\lambda}\right),\notag\\
&(c+j>0;~\lambda>0;~c\neq 0,-1,\dots;~j=0,\pm 1,\dots;~|w|+|z|<\lambda),\notag\\
\operatorname {J}_c^{jj}(c+j,c;\lambda,w,z)&=\int_0^\infty x^{c+j-1}e^{-\lambda x}{}_1F_1(c+j;c;wx){}_1F_1(c;c+ j;zx)dx\notag\\
&=\dfrac{\Gamma(c+j)}{(\lambda-w)^{j}(\lambda-z-w)^{c}}F_1\left(-j,j,c;c;\dfrac{w}{w-\lambda},\dfrac{w}{w+z-\lambda}\right),\notag\\
&(c+j>0;~\lambda>0;~c\neq 0,-1,\dots;~j=0,\pm 1,\dots;~|w|+|z|<\lambda),
\end{align}
where $F_1$ is the first Appell function (\cite{appell}, equation (1)). By mean of (\cite{slater}, formula (8.3.5))
\begin{align}\label{def6}
F_1(a;b,b';c;w,z)&=(1-w)^{-a}F_1\left(a,c-b-b',b';c;\dfrac{w}{w-1},\dfrac{z-w}{1-w}\right),
\end{align} 
it follows, for $j\geq p$, 
\begin{align}\label{integral7}
\operatorname {J}_c^{j(\pm p)}(b,b';\lambda,w,z)&=\int_0^\infty x^{c+j-1}e^{-\lambda x}{}_1F_1(b;c;wx){}_1F_1(b';c\pm p;zx)dx\notag\\
&=\dfrac{\Gamma(c+j)}{\lambda^{c+j-b-b'}(\lambda-w)^b(\lambda-z)^{b'}}\sum_{k=0}^{j\mp p}\dfrac{(-j\pm p)_k(b')_k}{(c\pm p)_k\,k!}\left(1-\dfrac{\lambda}{z}\right)^{-k}\notag\\
&\times \sum_{r=0}^{j+k}\dfrac{(b)_r\,(-j-k)_r}{(c)_r\,r!}\left(1-\dfrac{\lambda}{w}\right)^{-r}{}_2F_1\left(b+r,b'+k;c+r;\dfrac{w\,z}{(\lambda-z)(\lambda-w)}\right),\notag\\
&(c+j>0;~\lambda>0;~c,c\pm p\neq 0,-1,-2,\dots;~ |w|+|z|<\lambda),
\end{align}
where for $p=j= 0,1,2,\dots$
\begin{align}\label{integral8}
\operatorname {J}_c^{jj}(b,b';&\lambda,w,z)=\int_0^\infty x^{c+j-1}e^{-\lambda x}{}_1F_1(b;c;wx){}_1F_1(b';c+ j;zx)dx\notag\\
&=\dfrac{\Gamma(c+j)}{\lambda^{c+j-b-b'}(\lambda-w)^b(\lambda-z)^{b'}}
\sum_{r=0}^{j}\dfrac{(b)_r\,(-j)_r}{(c)_r\,r!}\left(1-\dfrac{\lambda}{w}\right)^{-r}{}_2F_1\left(b+r,b';c+r;\dfrac{w\,z}{(\lambda-z)(\lambda-w)}\right),\notag\\
&(c+j>0;~\lambda>0;~c,c+j\neq 0,-1,-2,\dots;~ |w|+|z|<\lambda).
\end{align}
Setting $b'=c+j$, equation \eqref{integral4} yield
\begin{align}\label{integral9}
\operatorname {J}_c^{j(\pm p)}(b,c+j;\lambda,w,z)&=\int_0^\infty x^{c+j-1}e^{-\lambda x}{}_1F_1(b;c;wx){}_1F_1(c+j;c\pm p;zx)dx\notag\\
&=\dfrac{\Gamma(c+j)}{(\lambda-z)^{c+j-b}(\lambda-z-w)^b}\sum_{k=0}^{j\mp p} \dfrac{(-j\pm p)_k(c+j)_k}{(c\pm p)_k~k!}\left(1-\dfrac{\lambda}{z}\right)^{-k}{}_2F_1\left(b,-j-k;c;\dfrac{w}{w+z-\lambda}\right),\notag\\
&(c+j>0;~\lambda>0;~c,c\pm p\neq 0,-1,-2,\dots;~p\geq 0;~j=0,\pm 1,\pm 2,\dots;~|w|+|z|<\lambda),
\end{align}
By means of the Kummer's first transformation
$
{}_1F_1(b;c;z)=e^z{}_1F_1(c-b;c;-z),
$
and the series representation
\begin{equation}\label{def10}
F_2(a; b,b';c,c'; x,y)= \sum_{m=0}^{\infty} {(a)_{m}\,(b)_m\over (c)_m\,m!}\, x^m\, {}_2F_1(a+m,b';c',y),
\end{equation}
it easily follows 
\begin{align}\label{integral11}
\operatorname {J}_c^{j(\pm p)}(c+j,b;\lambda,w,z)&=\int_0^\infty x^{c+j-1}e^{-\lambda x}{}_1F_1(c+j;c;wx){}_1F_1(b;c\pm p;zx)\,dx\notag\\
&=\dfrac{\Gamma(c+j)}{(\lambda-w)^{c+j}}\sum_{k=0}^j\dfrac{(-j)_k(c+j)_k}{(c)_k\,k!}\left(1-\dfrac{\lambda}{w}\right)^{-k}{}_2F_1\left(\begin{matrix}b,&c+j+k\\ 
~&c\pm p\end{matrix};\frac{z}{\lambda-w}\right),\notag\\
&(c+j>0;~\lambda>0;~c,c\pm p\neq0,-1,\dots;~|w|+|z|<\lambda).
\end{align}
By means of the identity (\cite{brychkov2008}, formula 5.14.3)
\begin{equation}\label{def12}
\sum_{k=0}^n{n\choose k} \dfrac{(-z)^k}{(b)_k}{}_1F_1(a;b+k;z)={}_1F_1(a-n;b;z),
\end{equation}
it follows that
\begin{align}\label{integral13}
\operatorname {J}_c^{j(\pm p)}(c-j,b;\lambda,w,z)&=\int_0^\infty x^{c+j-1}e^{-\lambda x}{}_1F_1(c-j;c;wx){}_1F_1(b;c\pm p;zx)\,dx\notag\\
&=\dfrac{\Gamma(c+j)}{\lambda^{c+j}}\sum_{k=0}^j\dfrac{(-j)_k(c+j)_k}{(c)_k\,k!}\left(\dfrac{w}{\lambda}\right)^kF_2\left(\begin{matrix}c+j+k,&c,&b\\
~&c+k&c\pm p\end{matrix};\frac{w}{\lambda},\frac{z}{\lambda}\right),\notag\\
&(c+j>0;~\lambda>0;~c,c\pm p\neq0,-1,\dots;~|w|+|z|<\lambda).
\end{align}
By means of the identity (\cite{brychkov2008}, formula 5.14.1)
\begin{equation}\label{def14}
\sum_{k=0}^n(-1)^k{n\choose k} \dfrac{(b-a)_k}{(b)_k}{}_1F_1(a;b+k;z)=\dfrac{(a)_n}{(b)_n}{}_1F_1(a+n;b+n;z),
\end{equation}
it follows that
\begin{align}\label{intgral15}
\operatorname {J}_{c+n}^{(j-n)(\pm p-n)}(b+n,b';\lambda,w,z)&=\int_0^\infty x^{c+j-1}e^{-\lambda x}{}_1F_1(b+n;c+n;wx){}_1F_1(b';c\pm p;z\,x)\,dx\notag\\
&=\dfrac{\Gamma(c+j)}{\lambda^{c+j}}\dfrac{(c)_n}{(b)_n}\sum_{k=0}^n\dfrac{(-n)_k(c-b)_k}{(c)_k\,k!}
F_2\left(\begin{matrix}c+j;&b,&b\rq{}\\ 
~&c+k,&c\pm p\end{matrix};\frac{w}{\lambda},\frac{z}{\lambda}\right),\notag\\
&(c+j>0;~\lambda>0;n=0,1,\dots;~c+n,c\pm p\neq0,-1,\dots;~|w|+|z|<\lambda).
\end{align}
 By means of the identity (\cite{brychkov2008}, formula 5.14.5)
\begin{equation}\label{def16}
{}_1F_1(b+n;c+n;w)=\dfrac{(c-1)_n(c)_n}{(b)_n(-w)^n}\sum_{k=0}^n\dfrac{(-n)_k(1-c)_k}{(2-c-n)_k\, k!}{}_1F_1(b,c-k,w)
\end{equation}
it follows that
\begin{align}\label{integral17}
\operatorname {J}_{c+n}^{(j-n)(\pm p-n)}&(b+n,b';\lambda,w,z)=\int_0^\infty x^{c+j-1}e^{-\lambda x}{}_1F_1(b+n;c+n;wx){}_1F_1(b';c\pm p;z\,x)\,dx\notag\\
&=\dfrac{(c-1)_n(c)_n}{(-w)^n\,(b)_n}\dfrac{\Gamma(c+j-n)}{\lambda^{c+j-n}}\sum_{k=0}^n\dfrac{(-n)_k(1-c)_k}{(2-c-n)_k\,k!}
F_2\left(\begin{matrix}c+j-n;&b,&b\rq{}\\ 
~&c-k,&c\pm p\end{matrix};\frac{w}{\lambda},\frac{z}{\lambda}\right),\notag\\
&(c\neq0,\pm 1;\dots,c+j>0;~\lambda>0;n=0,1,\dots;~c+n,c\pm p\neq0,-1,\dots;~|w|+|z|<\lambda).
\end{align}
By means of the identity (\cite{brychkov2008}, formula 5.14.6)
\begin{equation}\label{def18}
{}_1F_1(b+n;c;w)=\dfrac{(b-c+1)_n}{(b)_n}\sum_{k=0}^n(-1)^k{n\choose k}\dfrac{(1-c)_k}{(b-c+1)_k}{}_1F_1(b,c-k,w)
\end{equation}
it follows
\begin{align}\label{integral19}
\operatorname {J}_{c}^{j(\pm p)}(b+n,b';\lambda,w,z)&=\int_0^\infty x^{c+j-1}e^{-\lambda x}{}_1F_1(b+n;c;wx){}_1F_1(b';c\pm p;z\,x)\,dx\notag\\
&=\dfrac{(b-c+1)_n}{(b)_n}\dfrac{\Gamma(c+j)}{\lambda^{c+j}}\sum_{k=0}^n\dfrac{(-n)_k(1-c)_k}{(b-c+1)_k\,k!}
F_2\left(\begin{matrix}c+j;&b,&b\rq{}\\ 
~&c-k,&c\pm p\end{matrix};\frac{w}{\lambda},\frac{z}{\lambda}\right)\notag\\
&(c\neq 0,\pm 1,\dots,c+j>0;~\lambda>0;n=0,1,\dots;~c,c\pm p\neq0,-1,\dots;~|w|+|z|<\lambda).
\end{align}
By means of the identity (\cite{brychkov2008}, formula 5.14.7)
\begin{equation}\label{def20}
{}_1F_1(b-n;c-n;w)=\dfrac{(w)^n}{(1-c)_n}\sum_{k=0}^n{n\choose k}\dfrac{(1-c)_k}{w^k}{}_1F_1(b,c-k,w)
\end{equation}
it follows
\begin{align}\label{integral21}
\operatorname {J}_{c-n}^{(j+n)(\pm p+n)}&(b-n,b';\lambda,w,z)=\int_0^\infty x^{c+j-1}e^{-\lambda x}{}_1F_1(b-n;c-n;wx){}_1F_1(b';c\pm p;z\,x)\,dx\notag\\
&=\dfrac{w^n\Gamma(c+j+n)}{\lambda^{c+j+n}(1-c)_n}
\sum_{k=0}^n\dfrac{(-n)_k(1-c)_k}{k!\,(1-c-j-n)_k}
\left(\dfrac{\lambda}{w}\right)^k
F_2\left(\begin{matrix}c+j-k+n;&b,&b\rq{}\\ 
~&c-k,&c\pm p\end{matrix};\frac{w}{\lambda},\frac{z}{\lambda}\right)\notag\\
&(c\neq 0,\pm 1,\dots,c+j>0;~\lambda>0;n=0,1,\dots;~c-n,c\pm p\neq0,-1,\dots;~|w|+|z|<\lambda).
\end{align}
Since  (\cite{prudnikov}, formula 7.2.4.68)
\begin{equation}\label{def22}
F_2(a;b,b;c,c;z,-z)={}_4F_3\left(\dfrac{a}{2},\dfrac{a+1}{2},b,c-b;\dfrac{c}{2},\dfrac{c+2}{2},c;z^2\right).
\end{equation}
it follows that
\begin{align}
\operatorname {J}_c^{j0}(b,b;\lambda,w,-w)&=\int\limits_0^\infty x^{c+j-1}\,e^{-\lambda x}\,{}_1F_1(b;c;wx){}_1F_1(b;c;-wx)\,dx=\dfrac{\Gamma(c+j)}{\lambda^{c+j}}{}_4F_3\left(\begin{matrix}b,&c-b,&\dfrac{c+j}{2},&\dfrac{c+j+1}{2}\\ 
c,&\dfrac{c}{2},&\dfrac{c+1}{2}\end{matrix};\frac{w^2}{\lambda^2}\right),\notag\\
&(c+j>0;~\lambda>0;~c\neq 0,-1,-2,\dots,~|w|<\lambda),\label{integral23}\\
\operatorname {J}_c^{10}(b,b;\lambda;w,-w)&=\int\limits_0^\infty x^{c}\,e^{-\lambda x}\,{}_1F_1(b;c;w\,x)\,{}_1F_1(b;c;-w\,x)\,dx=\dfrac{\Gamma(c+1)}{\lambda^{c+1}}{}_3F_2\left(\begin{matrix}b,&c-b,&\dfrac{c}{2}+1\\ 
c,&\dfrac{c}{2},\end{matrix};\frac{w^2}{\lambda^2}\right),\notag\\
&(c>-1;~\lambda>0;~|w|<\lambda),\label{integral24}\\
\operatorname {J}_c^{10}\left(\dfrac{c}{2},\dfrac{c}{2};\lambda;w,-w\right)&=\int\limits_0^\infty x^{c}\,e^{-\lambda x}\,{}_1F_1\left(\dfrac{c}{2};c;w\,x\right)\,{}_1F_1\left(\dfrac{c}{2};c;-w\,x\right)\,dx=\dfrac{\Gamma(c+1)}{\lambda^{c+1}}{}_2F_1\left(\begin{matrix}\dfrac{c}{2},&\dfrac{c}{2}+1\\ \\
c,\end{matrix};\frac{w^2}{\lambda^2}\right),\notag\\
&(c>-1;~\lambda>0;~|w|<\lambda).\label{integral25}
\end{align}
On other hand, by means of (\cite{prudnikov}, formula 7.2.4.68)
\begin{equation}\label{def26}
F_2(a;b,c-b;c,c;z,z)=(1-z)^{-a}{}_4F_3\left(\dfrac{a}{2},~\dfrac{a+1}{2},~b,~c-b;~\dfrac{c}{2},~\dfrac{c+2}{2},~c;\dfrac{z^2}{(1-z)^2}\right),
\end{equation}
it follows that
\begin{align}
\operatorname {J}_c^{j0}(b,c-b;\lambda,z,z)&=\int\limits_0^\infty x^{c+j-1}e^{-\lambda x}{}_1F_1(b;c;z\,x){}_1F_1(c-b;c;z\,x)\,dx=\dfrac{\Gamma(c+j)}{(\lambda-z)^{c+j}}{}_4F_3\left(\begin{matrix}b,&c-b,&\dfrac{c+j}{2},&\dfrac{c+j+1}{2}\\ 
c,&\dfrac{c}{2},&\dfrac{c+1}{2}\end{matrix};\frac{z^2}{(\lambda-z)^2}\right),\notag\\
&(c+j>0;\lambda>0;~c\neq 0,-1,-2,\dots,~|w|<\lambda),\label{integral27}\\
\operatorname {J}_c^{10}(b,c-b;\lambda,z,z)&=\int\limits_0^\infty x^{c}\,e^{-\lambda x}\,{}_1F_1(b;c;z\,x)\,{}_1F_1(c-b;c;z\,x)\,dx=\dfrac{\Gamma(c+1)}{(\lambda-z)^{c+1}}{}_3F_2\left(\begin{matrix}b,&c-b,&\dfrac{c}{2}+1\\ 
c,&\dfrac{c}{2},\end{matrix};\frac{z^2}{(\lambda-z)^2}\right),\label{integral28}\\
&(c>-1;~\lambda>0;~|z|<|\lambda|),\notag\\
\operatorname {J}_c^{10}\left(\dfrac{c}{2},\dfrac{c}{2};\lambda,z,z\right)&=\int\limits_0^\infty x^{c}\,e^{-\lambda x}\left[{}_1F_1\left(\dfrac{c}{2};c;z\,x\right)\right]^2dx=\dfrac{\Gamma(c+1)}{(\lambda-z)^{c+1}}{}_2F_1\left(\begin{matrix}\dfrac{c}{2},&\dfrac{c}{2}+1\\ 
c,\end{matrix};\frac{z^2}{(\lambda-z)^2}\right),\notag\\
&(c>-1;~\lambda>0;~|z|<\lambda).\label{integral29}
\end{align}
By means of the identity (\cite{opps2009}, Theorem 3, formula 29)
\begin{align}\label{def30}
F_2(\sigma;\alpha_1,\alpha_2;\beta_1,\beta_2+n;w,z)&=  \frac{(\beta_2)_n}{(\beta_2-\alpha_2)_n}\sum\limits_{k=0}^n(-1)^k{n\choose k}{(\alpha_2)_k\over (\beta_2)_k} F_2(\sigma;\alpha_1,\alpha_2+k;\beta_1,\beta_2+k;w,z)\notag\\
&\left(|w|+|y|<1; n=0,1,2,\dots;\beta_1,\beta_2\neq
0,-1,-2,\dots;\beta_2>\alpha_2\right),
\end{align}
it easily follow 
\begin{align}\label{integral31}
\operatorname {J}_c^{jp}(b,b';\lambda;w,z)
&=\int_0^\infty x^{c+j-1}e^{-\lambda x}{}_1F_1(b;c;wx){}_1F_1(b';c+ p;zx)dx\notag\\
&=\frac{\Gamma(c+j)\,(c)_p}{\lambda^{c+j}(c-b')_p}\sum\limits_{k=0}^p\dfrac{(-p)_k\,(b')_k}{k!\,(c)_k} F_2\left(c+j;b,b'+k;c,c+k;\frac{w}{\lambda},\frac{z}{\lambda}\right),\notag\\
(c+j>0;p\geq 0;&\lambda>0;;~c \neq 0,-1,-2,\dots;~if~~c-b' (negative~integer), b'-c\geq p~|w|+|z|<\lambda),
\end{align} 
Further by means of (\cite{opps2009}, Theorem 3, formula 29)
\begin{align}\label{def32}
F_2\left(\begin{matrix}\sigma;&\alpha_1,&\alpha_2\\ 
~&\beta_1,&\beta_2-n\end{matrix};w,z\right)
&=\frac{1}{\left[\prod\limits_{i=0}^n(\beta_2-i)\right]}\sum\limits_{k=0}^n{n\choose k}\left[\prod\limits_{j=0}^{n-k}(\beta_2-j)\right]{(\sigma)_k(\alpha_2)_k\over (\beta_2)_k}z^k
F_2\left(\begin{matrix}\sigma+k;&\alpha_1,&\alpha_2+k\\ 
~&\beta_1,&\beta_2+k\end{matrix};w,z\right),\notag\\
 &\left(\beta_2\neq i, i=0,1,\dots,n;\beta_1,\beta_2-n\neq 0,-1,-2,\dots;n=0,1,2,\dots;|w|+|z|<1
 \right),
\end{align}
it follows
\begin{align}\label{integral33}
\operatorname {J}_c^{j(-p)}(b,b';\lambda;w,z)&
=\int_0^\infty x^{c+j-1}e^{-\lambda x}{}_1F_1(b;c;wx){}_1F_1(b';c-p;zx)dx\notag\\
&=\dfrac{\Gamma(c+j)}{\lambda^{c+j}}\sum\limits_{k=0}^p
\dfrac{(-p)_k\,(b')_k\,(c+j)_k}{(c)_k\,(c-p)_k\,k!}\left(-\dfrac{z}{\lambda}\right)^k
F_2\left(\begin{matrix}c+k+j;&b,&b'+k\\ 
~&c,&c+k\end{matrix};\frac{w}{\lambda},\frac{z}{\lambda}\right),\notag\\
&(c+j>0;~p\geq 0;\lambda>0;~c\neq 0,-1,-2,\dots;~c-p\neq 0,-1,-2,\dots;~|w|+|z|<\lambda)
\end{align}
whence
\begin{align}\label{integral34}
\operatorname {J}_c^{j(-p)}(b,c;\lambda;w,z)&
=\int_0^\infty x^{c+j-1}e^{-\lambda x}{}_1F_1(b;c;wx){}_1F_1(c;c-p;zx)dx\notag\\
&=\dfrac{\Gamma(c+j){(\lambda-z)^{b-c-j}}}{(\lambda-w-z)^b}
\sum\limits_{k=0}^p
\dfrac{(-p)_k(c+j)_k}{(c-p)_k\,k!}\left(\dfrac{z}{z-\lambda}\right)^k
{}_2F_1\left(\begin{matrix}-k-j,&b\\ 
~&c\end{matrix};\dfrac{w}{w+z-\lambda}\right)\notag\\
&(c+j>0;~p\geq 0;\lambda>0;~|w|+|z|<\lambda;~c\neq 0,-1,-2,\dots;~c-p\neq 0,-1,-2,\dots).
\end{align} 
The following identities are straightforward consequences of the pervious integrals:
\begin{align}\label{integral35}
\operatorname {J}_c^{j0}(c+j,0;\lambda,w,0)&=\int_0^\infty x^{c+j-1}\,e^{-\lambda\, x}{}_1F_1(c+j;~c;~w\,x)\,dx=\dfrac{\Gamma(c+j)}{(\lambda-w)^{c+j}}{}_2F_1\left(\begin{matrix}-j,&c+j\\ 
~&c\end{matrix};\frac{w}{w-\lambda}\right),\notag\\
&(c+j>0;~\lambda>0;~c\neq 0,-1,-2,\dots,~|w|<\lambda),
\end{align}
\begin{align}\label{integral36}
\operatorname {J}_c^{j(\pm p)}(0,b;\lambda,0,z)&=\int_0^\infty x^{c+j-1}\,e^{-\lambda x}\,{}_1F_1(b;c\pm p;z\,x)\, dx\notag\\
&=\dfrac{\Gamma(c+j)}{\lambda^{c+j}}{}_2F_1\left(\begin{matrix}c+j,&b\\ 
~&c\pm p\end{matrix};\frac{z}{\lambda}\right)=\dfrac{\Gamma(c+j)}{\lambda^{\pm p-b+ c}(\lambda-z)^{b\mp p+j}}{}_2F_1\left(\begin{matrix}\pm p-j,&c\pm p-b\\ 
~&c\pm p\end{matrix};\frac{z}{\lambda}\right),\notag\\
&(c+j>0;~\lambda>0;~c,c\pm p\neq 0,-1,-2,\dots;j=0,\pm 1,\dots;~p=0,1,2\dots;~|z|<\lambda),
\end{align} 
\begin{align}\label{integral37}
\operatorname {J}_c^{jj}(0,b;\lambda,0,z)&=\int_0^\infty x^{c+j-1}\,e^{-\lambda x}\,{}_1F_1(b;c+j;z\,x)\, dx=\dfrac{\Gamma(c+j)}{\lambda^{c-b+j}(\lambda-z)^b},\notag\\
&(c+j>0;\lambda>0;j=0,\pm 1,\dots;~|z|<\lambda),
\end{align} 

\begin{align}\label{integral38}
\operatorname {J}_c^{j0}(0,b;\lambda,0,z)&=\int_0^\infty x^{c+j-1}\,e^{-\lambda x}\,{}_1F_1(b;c;zx)\, dx=\dfrac{\Gamma(c+j)}{\lambda^{c+j-b}(\lambda-x)^b}
\left[1
+\dfrac{b\,z}{c\left(\lambda-z\right)}\sum_{k=1}^j{}_2F_1\left(\begin{matrix}-j+k,&b+1\\ \\
~&c+1\end{matrix};\dfrac{z}{z-\lambda}\right)
\right],\notag\\
&(c+j>0;~\lambda>0;~c\neq 0,-1,-2,\dots;~|z|<\lambda).
\end{align} 
\begin{align}\label{integral39}
\operatorname {J}_c^{00}(0,b;\lambda,0,z)&=\int_0^\infty x^{c-1}\,e^{-\lambda x}\,{}_1F_1(b;c;zx)\, dx=\dfrac{
\lambda^{b-c}\, \Gamma(c)}{(\lambda-z)^{b}},\qquad(c>0;~\lambda>0;~|z|<\lambda).
\end{align} 

\section{Gordon's integral and confluent hypergeometric polynomials}
\noindent In the case of $\alpha=-n$, the confluent hypergeometric function ${}_1F_1(\alpha;\beta;z)$ reduces to $n$-degree polynomial in $z$, namely  ${}_1F_1(-n;\beta;z)=\sum_{k=0}^n{(-n)_k\,z^k}/((\beta)_k\,k!),~n=0,1,\dots.
$.
 Thus,
\begin{align}\label{integral40}
\operatorname {J}_c^{j(\pm p)}(b,-n;\lambda,w,z)&=\int_0^\infty x^{c+j-1}e^{-\lambda x}{}_1F_1(b;c;wx){}_1F_1(-n;c\pm p;zx)dx\notag\\
&=\dfrac{\Gamma(c+j)}{\lambda^{c+j-b}(\lambda-w)^b}
\sum_{k=0}^n\dfrac{(-n)_k(c+j)_k}{(c\pm p)_k\,k!}\left(\dfrac{z}{\lambda}\right)^k{}_2F_1(-j-k,b;c;\dfrac{w}{w-\lambda}),\notag\\
&(c+j>0;~\lambda>0;~c,~c\pm p\neq 0,-1,\dots;p=0,1,\dots;~|w|<|\lambda|),
\end{align}
where a direct differentiation of both sides with respect to $z$ yields
\begin{align}\label{integral41}
\operatorname {J}_c^{(j+m)(\pm p+m)}(b,m-n;\lambda,w,z)&=\int_0^\infty x^{c+j+m-1}e^{-\lambda x}{}_1F_1(b;c;wx){}_1F_1(m-n;c\pm p+m;zx)\,dx\notag\\
&=\dfrac{(-1)^m\Gamma(c+j)(c\pm p)_m}{(-n)_mz^m\lambda^{c+j-b}(\lambda-w)^b}
\sum_{k=m}^n\dfrac{(-k)_m(-n)_k(c+j)_k}{(c\pm p)_k\,k!}\left(\dfrac{z}{\lambda}\right)^k{}_2F_1(-j-k,b;c;\dfrac{w}{w-\lambda}),\notag\\
&(m\leq n; c+j+m>0;~\lambda>0;~c,c,~c\pm p+m\neq 0,-1,\dots;~|w|<|\lambda|).
\end{align}
Further, setting $b=-m$ in equation \eqref{integral40} implies
\begin{align}\label{integral42}
\operatorname {J}_c^{jp}(-m,-n;\lambda;w,z)&=\int_0^\infty x^{c+j-1}e^{-\lambda\, x}{}_1F_1(-m;c;w\,x)\,{}_1F_1(-n;c\pm p;z\,x)\,dx\notag\\
&=\dfrac{\Gamma(c+j)}{\lambda^{c+j}}\sum_{k=0}^m\dfrac{(c+j)_k(-m)_k}{(c)_k\,k!}\left(\dfrac{w}{\lambda}\right)^k{}_2F_1(-n,c+j+k;c\pm p;\dfrac{z}{\lambda})\notag\\
&\equiv\dfrac{\Gamma(c+j)}{\lambda^{c+j}}\sum_{k=0}^n\dfrac{(c+j)_k(-n)_k}{(c\pm p)_k\,k!}\left(\dfrac{z}{\lambda}\right)^k{}_2F_1(-m,c+j+k;c;\dfrac{w}{\lambda}),\notag\\
&(c+j>0,~\lambda>0;~c,c\pm p\neq 0,-1,\dots;j=0,\pm 1,\dots; n,m=0,1,\dots),
\end{align}
where a direct differentiation of both sides with respect to $w$ yields
\begin{align}\label{integral43}
\int_0^\infty x^{c+j+l-1}&e^{-\lambda\, x}{}_1F_1(l-m;c+l;w\,x)\,{}_1F_1(-n;c\pm p;z\,x)\,dx\notag\\
&=\dfrac{(-1)^l\Gamma(c+j)(c)_l}{\lambda^{c+j}w^l(-m)_l}\sum_{k=l}^m\dfrac{(c+j)_k(-m)_k(-k)_l}{(c)_k\,k!}\left(\dfrac{w}{\lambda}\right)^k{}_2F_1(-n,c+j+k;c\pm p;\dfrac{z}{\lambda}),\notag\\
&(l\leq m;c+j+l>0,~\lambda>0;~c+l,c\pm p\neq 0,-1,\dots;j=0,\pm 1,\dots; n,m=0,1,\dots),
\end{align}
with a further differentiation of both sides with respect to $z$ yields
\begin{align}\label{integral44}
\int_0^\infty& x^{c+j+k+s-1}e^{-\lambda\, x}{}_1F_1(l-m;c+l;w\,x)\,{}_1F_1(s-n;c\pm p+s;z\,x)\,dx\notag\\
&=\dfrac{(-1)^l\Gamma(c+j)(c)_l}{\lambda^{c+j+s}w^l(-m)_l}{}\sum_{k=l}^m\dfrac{(c+j+k)_s(c+j)_k(-m)_k(-k)_l}{(c)_k\,k!}\left(\dfrac{w}{\lambda}\right)^k{}_2F_1(s-n,c+j+k+s;c\pm p+s;\dfrac{z}{\lambda})\notag\\
&(s\leq n;l\leq m;c+j+l+s>0,~\lambda>0;~c+l,c\pm p+s\neq 0,-1,\dots;s,l,j=0,\pm 1,\dots; n,m=0,1,\dots).
\end{align}
If $z=\lambda$, equation \eqref{integral40} reads
\begin{align}\label{integral45}
\operatorname {J}_c^{j(\pm p)}(-m,-n;\lambda;w,\lambda)&=\int_0^\infty x^{c+j-1}e^{-\lambda\, x}{}_1F_1(-m;c;w\,x)\,{}_1F_1(-n;c\pm p;\lambda\,x)\,dx\notag\\
&=\dfrac{\Gamma(c+j)(\pm p-j)_n}{\lambda^{c+j}(\pm p+c)_n}{}_3F_2(-m,c+j,1+j\mp p;c,1+j-n\mp p;\dfrac{w}{\lambda}),\notag\\
(c+j>0,~\lambda>0;~c,c\pm p&\neq 0,-1,\dots;j=0,\pm 1,\dots; n,m=0,1,\dots;1+j-n\pm p\neq 0,-1,\dots),
\end{align}
and  if $w=\lambda$
equation \eqref{integral40} reads
\begin{align}\label{integral46}
\operatorname {J}_c^{j(\pm p)}(-m,-n;\lambda;\lambda,z)&=\int_0^\infty x^{c+j-1}e^{-\lambda\, x}{}_1F_1(-m;c;\lambda\,x)\,{}_1F_1(-n;c\pm p;z\,x)\,dx\notag\\
&=\dfrac{\Gamma(c+j)(-j)_m}{\lambda^{c+j}(c)_m}{}_3F_2(-n,c+j,1+j;c\pm p,1+j-m;\dfrac{z}{\lambda}),\notag\\
(c+j>0,~\lambda>0;~c,c\pm p&\neq 0,-1,\dots;j=0,\pm 1,\dots; n,m=0,1,\dots;1+j-m\neq 0,-1,\dots).
\end{align}
Further if $w=\lambda$, equation \eqref{integral45} reads
\begin{align}\label{integral47}
\operatorname {J}_c^{jp}(-m,-n;\lambda;\lambda,\lambda)&=\int_0^\infty x^{c+j-1}e^{-\lambda\, x}{}_1F_1(-m;c;\lambda\,x)\,{}_1F_1(-n;c\pm p;\lambda\,x)\,dx\notag\\
&=\dfrac{\Gamma(c+j)\,(\pm p-j)_n}{\lambda^{c+j}\,(\pm p+c)_n}\,{}_3F_2(-m,c+j,1+j\mp p;c,1+j-n\mp p;1)\notag\\
&=\dfrac{\Gamma(c+j)(-j)_m}{\lambda^{c+j}(c)_m}{}_3F_2(-n,c+j,1+j;c\pm p,1+j-m;1),\notag\\
(c+j>0,~\lambda>0;~c,c\pm p&\neq 0,-1,\dots;j=0,\pm 1,\dots; n,m=0,1,\dots;1+j-n\pm p,1+j-m\neq 0,-1,\dots),
\end{align}
whenece
\begin{align}\label{integral48}
\operatorname {J}_c^{j0}(-m,-n;\lambda;\lambda,\lambda)&=\int_0^\infty x^{c+j-1}e^{-\lambda\, x}{}_1F_1(-m;c;\lambda\,x)\,{}_1F_1(-n;c;\lambda\,x)\,dx\notag\\
&=\dfrac{\Gamma(c+j)\,(-j)_n}{\lambda^{c+j}\,(c)_n}\,{}_3F_2(-m,c+j,1+j;c,1+j-n;1),\notag\\
&\equiv\dfrac{\Gamma(c+j)\,(-j)_m}{\lambda^{c+j}\,(c)_m}\,{}_3F_2(-n,c+j,1+j;c,1+j-m;1),\notag\\
(c+j>0,~\lambda>0;~c\neq 0,-1,&\dots;j=0,\pm 1,\dots; n,m=0,1,\dots;1+j-n,1+j-m\neq 0,-1,\dots).
\end{align}
If $j=0$, equation \eqref{integral47}
\begin{align}\label{integral49}
\operatorname {J}_c^{0p_\pm}(-n,-m;\lambda,\lambda,\lambda)&=\int_0^\infty x^{c-1}e^{-\lambda x}{}_1F_1(-n;c;\lambda\,x){}_1F_1(-m;c\pm p;\lambda x)\,dx=\dfrac{\Gamma(c)}{\lambda^{c}}\dfrac{m!}{(m-n)!} \dfrac{ (\pm p)_{m-n}}{(c\pm p)_m},\notag\\
(m\geq n;~c>0,~\lambda>0;&~c\pm p\neq 0,-1,\dots;j=0,\pm 1,\dots; n,m=0,1,\dots).
\end{align} 
If $m=n$, equation \eqref{integral47}
\begin{align}\label{integral50}
\operatorname {J}_c^{j0}(-n,-n;\lambda;\lambda,\lambda)&=\int_0^\infty x^{c+j-1}e^{-\lambda\, x}\left[{}_1F_1(-n;c;\lambda\,x)\,\right]^2\,dx=\dfrac{\Gamma(c+j)\,(-j)_n}{\lambda^{c+j}\,(c)_n}\,{}_3F_2(-n,c+j,1+j;c,1+j-n;1),\notag\\
(c+j>0,~\lambda>0;~c&\neq 0,-1,\dots;j=0,\pm 1,\dots; n=0,1,\dots;1+j-n\neq 0,-1,\dots).
\end{align}
The condition $1+j-n\neq 0,-1,-2,\dots$ in \eqref{integral50} can be softened using the identity
\begin{align}\label{eq51}
(-j)_n\,{}_3F_2(-n,c+j,1+j;c,1+j-n;1)=n!\,{}_3F_2(-n,-j,j+1;c,1;1),
\end{align}
to yield
\begin{align}\label{integral52}
\operatorname {J}_c^{j0}(-n,-n;\lambda;\lambda,\lambda)&=\int_0^\infty x^{c+j-1}e^{-\lambda\, x}\left[{}_1F_1(-n;c;\lambda\,x)\right]^2\,dx=\dfrac{\Gamma(c+j)\,n!}{\lambda^{c+j}\,(c)_n}\,{}_3F_2(-n,-j,1+j;c,1;1),\notag\\
&(c+j>0;~\lambda>0;~c\neq 0,-1,-2,\dots;j=0,\pm 1,\pm 2,\dots; n=0,1,2,\dots).
\end{align}
and thus
\begin{align}\label{integral53}
\operatorname {J}_c^{10}(-n,-n;\lambda;\lambda,\lambda)&=\int_0^\infty x^{c}\,e^{-\lambda\, x}\left[{}_1F_1(-n;c;\lambda\,x)\right]^2\,dx=\dfrac{\Gamma(c)\,n!}{\lambda^{c+1}\,(c)_n}\left(c+2n\right),~~(c>0;~\lambda>0).
\end{align}
From the symmetric property $j\longleftrightarrow -j-1$ of ${}_3F_2(-n,-j,1+j;c,1;1)$, it also follow
\begin{align}\label{integral54}
\operatorname {J}_c^{(-j-1)0}(-n,-n;\lambda;\lambda,\lambda)&=\int_0^\infty x^{c-j-2}e^{-\lambda\, x}\left[{}_1F_1(-n;c;\lambda\,x)\right]^2\,dx=\dfrac{\Gamma(c-j-1)\,n!}{\lambda^{c-j-1}\,(c)_n}\,{}_3F_2(-n,-j,1+j;c,1;1),\notag\\
&(c-j-2>0;~\lambda>0;~c\neq 0,-1,-2,\dots;j=0,\pm 1,\pm 2,\dots; n=0,1,2,\dots)
\end{align}
whence
\begin{align}\label{integral55}
\operatorname {J}_c^{(-j-1)0}(-n,-n;\lambda;\lambda,\lambda)
&=\dfrac{\Gamma(c-j-1)\lambda^{2j+1}}{\Gamma(c+j)}\operatorname {J}_c^{j0}(-n,-n;\lambda;\lambda,\lambda),
\end{align}
for example
\begin{align}\label{integral56}
\operatorname {J}_c^{(-2)0}(-n,-n;\lambda;\lambda,\lambda)&=\dfrac{\Gamma(c-2)\lambda^{3}}{\Gamma(c+1)}\operatorname {J}_c^{10}(-n,-n;\lambda;\lambda,\lambda)=\dfrac{\Gamma(c-2)n!}{c\lambda^{c-2}\,(c)_n}\left(c+2n\right).
\end{align}
Further recurrence relations of this type are developed in the appendix. Note, from equations \eqref{integral42} and \eqref{integral45}, it follows
\begin{align}\label{eq57}
\sum_{k=0}^n\dfrac{(-n)_k\,(c+j)_k}{(c\pm p)_k\,k!}{}_2F_1\left(-m,c+j+k;c;\dfrac{w}{\lambda}\right)=\dfrac{(\pm p-j)_n}{(\pm p+c)_n}\,{}_3F_2\left(-m,c+j,1+j\mp p;c,1+j-n\mp p;\dfrac{w}{\lambda}\right).
\end{align}
An important class of W. Gordon's integral occur in the case of $w=k_1,z=k_2$, and $\lambda=(k_1+k_2)/2$, namely,
\begin{align}\label{integral58}
\operatorname {J}_c^{j(\pm p)}\left(-n,-m;\dfrac{k_1+k_2}{2},k_1,k_2\right)&=\int_0^\infty x^{c+j-1}e^{-(k_1+k_2) x/2}{}_1F_1(-n;c;k_1\, x){}_1F_1(-m;c\pm p;k_2\, x)dx\notag\\
&=\dfrac{2^{c+j}\Gamma(c+j)}{(k_1+k_2)^{c+j}}\,F_2\left(\begin{matrix}c+j;&-n,&-m\\ 
~&c,&c\pm p\end{matrix};\frac{2k_1}{k_1+k_2},\frac{2k_2}{k_1+k_2}\right),\notag\\
&(k_1+k_2>0; c+j>0;c,c\pm p\neq 0,\pm 1,\dots)
\end{align} 
equivalently,
\begin{align}\label{integral59}
\operatorname {J}_c^{j(\pm p)}&\left(-n,-m;\dfrac{k_1+k_2}{2},k_1,k_2\right)=\int_0^\infty x^{c+j-1}e^{-(k_1+k_2) x/2}{}_1F_1(-n;c;k_1\, x){}_1F_1(-m;c\pm p;k_2\, x)dx\notag\\
&=\left\{ \begin{array}{ll}
\dfrac{\Gamma(c+j)(\pm p-j)_m}{k_1^{c+j}(\pm p+c)_m}{}_3F_2(-n,c+j,1+j\mp p;c,1+j-m\mp p;1),&\mbox{if $k_1= k_2,$} \\ \\
\dfrac{2^{c+j}\Gamma(c+j)}{(k_1+k_2)^{c+j}}\left(\dfrac{k_1-k_2}{k_1+k_2}\right)^m
\sum_{i=0}^{\min\{j\mp p,m\}}  \dfrac{(-j\pm p)_i(-m)_i}{(c\pm p)_i~i!}\left(\dfrac{2k_2}{k_2-k_1}\right)^{i}\\
\times F_1\left(-n,c+j+m,i-m;c;\dfrac{2k_1}{k_1+k_2},\dfrac{2k_1}{k_1-k_2}\right), &\mbox{if $k_1\neq k_2$.}
       \end{array} \right.
\end{align} 
and further equivalent to
\begin{align}\label{integral60}
\operatorname {J}_c^{j(\pm p)}&\left(-n,-m;\dfrac{k_1+k_2}{2},k_1,k_2\right)=\int_0^\infty x^{c+j-1}e^{-(k_1+k_2) x/2}{}_1F_1(-n;c;k_1\, x){}_1F_1(-m;c\pm p;k_2\, x)dx\notag\\
&=\left\{ \begin{array}{l}
\dfrac{\Gamma(c+j)(\pm p-j)_m}{k_1^{c+j}(\pm p+c)_m}{}_3F_2(-n,c+j,1+j\mp p;c,1+j-m\mp p;1),\mbox{if $k_1= k_2,$} \\ \\
\dfrac{(-1)^n2^{c+j}\Gamma(c+j)}{(k_1+k_2)^{c+j}}\left(\dfrac{k_1-k_2}{k_1+k_2}\right)^{m+n}
\sum_{i=0}^{\min\{j\mp p,m\}}  \dfrac{(-j\pm p)_i(-m)_i}{(c\pm p)_i~i!}\left(\dfrac{2k_2}{k_2-k_1}\right)^{i}\\
\times\sum\limits_{r=0}^{j+i}\dfrac{(-n)_r(-j-i)_r}{(c)_r\,r!}\left(\dfrac{2k_2}{k_2-k_1}\right)^r {}_2F_1\left(r-n,i-m;c+r;\dfrac{-4k_1k_2}{(k_1-k_2)^2}\right), \mbox{\qquad\qquad if $k_1\neq k_2$.}
       \end{array} \right.
\end{align} 
In particular
\begin{align}\label{integral61}
\operatorname {J}_c^{j0}&\left(-n,-m;\dfrac{k_1+k_2}{2},k_1,k_2\right)=\int_0^\infty x^{c+j-1}e^{-(k_1+k_2) x/2}{}_1F_1(-n;c;k_1\, x){}_1F_1(-m;c;k_2\, x)dx\notag\\
&=\left\{ \begin{array}{l}
\dfrac{\Gamma(c+j)(-j)_m}{k_1^{c+j}(\pm p+c)_m}{}_3F_2(-n,c+j,1+j;c,1+j-m;1),, \qquad\qquad (k_1= k_2, m\leq j),\\ \\
\dfrac{(-1)^n2^{c+j}\Gamma(c+j)}{(k_1+k_2)^{c+j}}\left(\dfrac{k_1-k_2}{k_1+k_2}\right)^{m+n}
\sum_{i=0}^m  \dfrac{(-j)_i(-m)_i}{(c)_i~i!}\left(\dfrac{2k_2}{k_2-k_1}\right)^{i}\\
\times\sum\limits_{r=0}^{j+i}\dfrac{(-n)_r(-j-i)_r}{(c)_r\,r!}\left(\dfrac{2k_2}{k_2-k_1}\right)^r {}_2F_1\left(r-n,i-m;c+r;\dfrac{-4k_1k_2}{(k_1-k_2)^2}\right), \quad\quad (k_1\neq k_2, m\leq j),
       \end{array} \right.
\end{align} 
and
\begin{align}\label{integral62}
\operatorname {J}_c^{00}\left(-m,-n;\dfrac{k_1+k_2}{2};k_1,k_2\right)&=\int_0^\infty x^{c-1}e^{-(k_1+k_2)x/2}{}_1F_1(-m;c;k_1\,x)\,{}_1F_1(-n;c;k_2\,x)\,dx\notag\\
&=\dfrac{2^c\Gamma(c)}{(k_1+k_2)^{c}}\sum_{k=0}^m\dfrac{(-m)_k}{k!}\left(\dfrac{2k_1}{k_1+k_2}\right)^k{}_2F_1(-n,c+k;c;\dfrac{2k_2}{k_1+k_2}),\notag\\
&(c>0,~\lambda>0;~c\neq 0,-1,-2,\dots; n,m=0,1,2,\dots),
\end{align}
where, generally, 
\begin{align}\label{integral63}
\operatorname {J}_c^{00}(-m,-n;\lambda;w,z)&=\int_0^\infty x^{c-1}e^{-\lambda\, x}{}_1F_1(-m;c;w\,x)\,{}_1F_1(-n;c;z\,x)\,dx\notag\\
&=\dfrac{\Gamma(c)}{\lambda^{c}}\sum_{k=0}^m\dfrac{(-m)_k}{k!}\left(\dfrac{w}{\lambda}\right)^k{}_2F_1(-n,c+k;c;\dfrac{z}{\lambda}),\notag\\
&(c>0,~\lambda>0;~c\neq 0,-1,-2,\dots; n,m=0,1,2,\dots),
\end{align}
from which the classical orthogonality property of the confluent hypergeometric functions follows, namely,
\begin{align}\label{integral64}
\operatorname {J}_c^{00}(-m,-n;\lambda;\lambda,\lambda)&=\int_0^\infty x^{c-1}e^{-\lambda\, x}{}_1F_1(-m;c;\lambda\,x)\,{}_1F_1(-n;c;\lambda\,x)\,dx=\dfrac{\Gamma(c)\,n!}{\lambda^c\,(c)_n}\delta_{nm},\notag\\
&(c>0,~\lambda>0;~c\neq 0,-1,-2,\dots; \delta_{nm}=0~if~ n\neq m,\delta_{nm}=1~if~ n= m),
\end{align}
using
$\sum_{k=0}^m{(-m)_k(-k)_n}/{k!}=n!\delta_{nm}$.
The same conclusion also follows from equation \eqref{integral48} using the fact that
$$\lim_{j\rightarrow 0}\, (-j)_m\,{}_3F_2(-n,c+j,1+j;c,1+j-m;1)=n!\,\delta_{nm}.
$$
If in equation \eqref{integral1}, $b=b'=-n$ and $p=0$, it follows 
\begin{align}\label{integral65}
\operatorname {J}_c^{j0}(-n,-n;\lambda;w,z)&=\int_0^\infty x^{c+j-1}e^{-\lambda\, x}{}_1F_1(-n;c;w\,x)\,{}_1F_1(-n;c;z\,x)\,dx\notag\\
&=\dfrac{n!\,\Gamma(c+j)}{\lambda^{c+j}(c)_n}\sum_{k=0}^n\dfrac{(c+j)_k(-n)_k}{(c)_k\,k!}\left(\dfrac{z}{\lambda}\right)^kP_n^{(c-1,j+k-n)}\left(1-\dfrac{2w}{\lambda}\right),
\end{align}
where $P_n^{(\alpha,\beta)}(z)$ is the Jacobi polynomial of order $\alpha$, $\beta$ and degree $n$ in $z$. 
The relation
$ P_n^{(a,b)}(-1)={(-1)^n}(b+1)_n/{n!}$ reduce the equation \eqref{integral65} to
\begin{align}\label{integral66}
\operatorname {J}_c^{j0}(-n,-n;\lambda;\lambda,z)&=\int_0^\infty x^{c+j-1}e^{-\lambda\, x}{}_1F_1(-n;c;\lambda\,x)\,{}_1F_1(-n;c;z\,x)\,dx\notag\\
&=\dfrac{\Gamma(c+j)(-j)_n}{\lambda^{c+j}\,(c)_n}\,{}_3F_2\left(-n,c+j,1+j;c,1+j-n;\dfrac{z}{\lambda}\right),\notag\\
&(c+j>0,~\lambda>0;~c\neq 0,-1,\dots; n=0,1,\dots;1+j-n\neq 0,-1,\dots),
\end{align}
as expected.
From equation \eqref{integral66}, it follows
\begin{align}\label{integral67}
\operatorname {J}_c^{n0}(-n,-n;\lambda;\lambda,z)=&\int_0^\infty x^{c+n-1}e^{-\lambda\, x}{}_1F_1(-n;c;\lambda\,x)\,{}_1F_1(-n;c;z\,x)\,dx=\dfrac{(-1)^n\,\Gamma(c)\,n!}{\lambda^{c+n}}\,{}_3F_2\left(-n,c+n,1+n;c,1;\dfrac{z}{\lambda}\right),\notag\\
&(c+n>0,~\lambda>0;~c\neq 0,-1,\dots).
\end{align}
For $n\geq m$
\begin{align}\label{integral68}
\operatorname {J}_c^{(n-m)(\pm p)}(-n,-m;\lambda,\lambda,z)&=\int_0^\infty x^{c+n-m-1}e^{-\lambda x}{}_1F_1(-n;c;\lambda\,x){}_1F_1(-m;c\pm p;z\,x)\,dx\notag\\
&=\dfrac{(-1)^{m+n}\,\Gamma(c)\,n!}{\lambda^{c+n-m}(c\pm p)_m}\left(\dfrac{z}{\lambda}\right)^m,\quad (c+n-m>0;\lambda>0;c,c\pm p\neq 0,-1,\dots;p\geq 0).
\end{align} 
The following integral follows immediately
\begin{align}\label{integral69}
\operatorname {J}_c^{j(\pm p)}(0,-n;\lambda,0,z)&=\int_0^\infty x^{c+j-1}e^{-\lambda x}{}_1F_1(-n;c\pm p;z\,x)dx=\dfrac{\Gamma(c+j)}{\lambda^{c+j}}
{}_2F_1\left(-n,c+j;c\pm p;\dfrac{z}{\lambda}\right),\notag\\
&(c+j>0;~\lambda>0;~c\pm p\neq 0,-1,\dots;p=0,1,\dots)
\end{align}
whence 
\begin{align}\label{integral70}
\operatorname {J}_c^{j(\pm p)}(0,-n;\lambda,0,\lambda)&=\int_0^\infty x^{c+j-1}e^{-\lambda x}{}_1F_1(-n;c\pm p;\lambda x)dx=\left\{ \begin{array}{ll}
\dfrac{\Gamma(c+j)}{\lambda^{c+j}}\dfrac{(\pm p-j)_n}{(c\pm p)_n}, &\mbox{ if $j\mp p\geq n,$} \\ 
0 ,&\mbox{ if $j\mp p< n,$}
       \end{array} \right.\notag\\
&(c+j>0;~\lambda>0;~c\pm p\neq 0,-1,\dots;p=0,1,\dots),
\end{align}
and if $p=j$,
\begin{align}\label{integral71}
\operatorname {J}_c^{jj}(0,-n;\lambda,0,z)&=\int_0^\infty x^{c+j-1}e^{-\lambda x}{}_1F_1(-n;c+j;z\, x)dx
=\dfrac{\Gamma(c+j)}{\lambda^{c+j}}\left(1-\dfrac{z}{\lambda}\right)^n,\quad( c+j>0;\lambda>0).
\end{align}
and
\begin{align}\label{integral72}
\operatorname {J}_c^{j0}(0,-n;\lambda,0,\lambda)&=\int_0^\infty x^{c+n-1}e^{-\lambda x}{}_1F_1(-n;c;\lambda\, x)dx=\dfrac{(-1)^n\,n!\,\Gamma(c)}{\lambda^{c+n}},\quad (c+n>0;~\lambda>0; ~c\neq 0,-1,\dots).
\end{align}

\section{Goron\rq{}s Integral and special functions}
\noindent The generalized Laguerre polynomials are defined, for integer $n$, in terms of confluent hypergeometric functions by 
\begin{align}\label{def73}
L_n^\lambda(z)=\dfrac{(\lambda+1)_n}{n!}{}_1F_1(-n;\lambda+1;z),
\end{align}
thus, 
\begin{align}\label{integral74}
\int_0^\infty x^{c+j-1}&e^{-\lambda x}
L_n^{c\pm p -1}(z\,x) {}_1F_1(b;c;w\,x)dx=\dfrac{\Gamma(c+j)(c\pm p)_n}{n!\lambda^{c+j}}
\sum_{k=0}^n\dfrac{(-n)_k(c+j)_k}{(c\pm p)_k\,k!}\left(\dfrac{z}{\lambda}\right)^k{}_2F_1\left(c+j+k,b;c;\dfrac{w}{\lambda}\right)\notag\\
&(c+j>0;~\lambda>0;~c,~c\pm p\neq 0,-1,\dots;p=0,1,\dots;~|w|<|\lambda|),
\end{align}
whence, if $b=0$,
\begin{align}\label{integral75}
\int_0^\infty x^{c+j-1}&e^{-\lambda x}
L_n^{c\pm p -1}(z\,x)dx=\left\{ \begin{array}{ll}
\dfrac{\Gamma(c+j)(c\pm p)_n}{n!\,\lambda^{c+j}}
{}_2F_1\left(-n,c+j;c\pm p;\dfrac{z}{\lambda}\right), &\mbox{ if $z\neq \lambda,~c\pm p\neq 0,-1,\dots,$} \\ 
\dfrac{\Gamma(c+j)(\pm p-j)_n}{\lambda^{c+j}\, n!} ,&\mbox{ if $z=\lambda ,~j\mp p\geq n,$}
       \end{array} \right.\notag\\
&(c+j>0;~\lambda>0;~p=0,1,\dots).
\end{align}

\noindent From equation \eqref{integral75}, it follows  
\begin{align}\label{integral76}
\int_0^\infty x^{c}\,e^{-\lambda x}\,
L_n^{c }(\lambda x)\,dx&= \left\{ \begin{array}{ll}
\dfrac{\Gamma(c+1)}{\lambda^{1+c}}, &\mbox{ if $n=0$}, \\ \\
0, &\mbox{ if $n\geq 1;c>-1,\lambda>0,~ n=0,1,\dots$},
       \end{array} \right.\notag\\ \notag\\
\int_0^\infty x^{c+j}\,e^{-\lambda x}\,
L_n^{c }(\lambda x)\,dx&= \left\{ \begin{array}{ll}
\dfrac{\Gamma(c+j+1)\, (-j)_n}{\lambda^{j+c+1}\, n!}, &\mbox{ if $n<j$}, \\ \\
(-1)^n\dfrac{\Gamma(c+n+1)}{\lambda^{n+c+1}} ,&\mbox{ if $n=j$}, \\ \\
0, &\mbox{ if $n> j; c+j>-1, ~\lambda>0,~n=0,1,\dots,$}
       \end{array} \right.\notag\\ \notag\\
\int_0^\infty x^{c}\,e^{-\lambda x}\,
L_n^{c- p }(\lambda x)\,dx&= \left\{ \begin{array}{ll}
\dfrac{\Gamma(c+1)(-p)_n}{\lambda^{c+1}\, n!},&\mbox{ if $n<p$},\\ \\
(-1)^n\dfrac{\Gamma(c+1)}{\lambda^{c+1}}, &\mbox{ if $n=p$}, \\ \\
0, &\mbox{ if $n> p;c>-1,\lambda>0,c-p\geq 0$.}
       \end{array} \right.
\end{align}
Since
$$\dfrac{d^m}{dz^m}L_n^\lambda(az)=(-a)^mL_{n-m}^{\lambda+m}(az)$$
it easily follows, for $n\geq m$ and $m=0,1,2,\dots$, that
\begin{align}\label{integral77}
\int_0^\infty x^{c+j+m-1}&e^{-\lambda x}
L_{n-m}^{c\pm p+m-1}(zx){}_1F_1(b;c;wx)dx\notag\\
&=\dfrac{(c\pm p)_n\Gamma(c+j)}{n!\,z^m\,\lambda^{c+j}}
\sum_{k=m}^n\dfrac{(-k)_m(-n)_k(c+j)_k}{(c\pm p)_k\,k!}\left(\dfrac{z}{\lambda}\right)^k{}_2F_1(c+j+k,b;c;\dfrac{w}{\lambda}).
\end{align}
and, for $\mu=0,1,2,\dots, m\leq n$,
\begin{align}\label{integral78}
\int_0^\infty &x^{c+j+m+\mu-1}e^{-\lambda x}
L_{n-m}^{c\pm p+m-1}(zx){}_1F_1(b+\mu;c+\mu;w\,x)dx\notag\\
&=\dfrac{\Gamma(c+j)(c+p)_n}{n!\,z^m\,\lambda^{c+j+\mu}}
\sum_{k=0}^n\dfrac{(-k)_m(-n)_k(c+j)_k(c+j+k)_\mu}{(c\pm p)_k\,k!}\left(\dfrac{z}{\lambda}\right)^k{}_2F_1(c+j+k+\mu,b+\mu;c+\mu;\dfrac{w}{\lambda}).
\end{align}
On other hand,
\begin{align}\label{integral79}
\int_0^\infty x^{c+j-1}e^{-\lambda x}
L_n^{c\pm p -1}(zx)L_m^{c-1}(w\, x)dx&=\dfrac{(c)_m(c\pm p)_n\Gamma(c+j)}{m!\, n!\,\lambda^{c+j}}
\sum_{k=0}^n\dfrac{(-n)_k(c+j)_k}{(c\pm p)_k\,k!}\left(\dfrac{z}{\lambda}\right)^k{}_2F_1\left(c+j+k,-m;c;\dfrac{w}{\lambda}\right),\notag\\
&(c+j>0;~\lambda>0;~c,~c\pm p\neq 0,-1,\dots;p=0,1,\dots).
\end{align}
and by direct differentiation $s$-times, with respect to $w$, of both sides
\begin{align}\label{integral80}
\int_0^\infty& x^{c+j+s-1}e^{-\lambda x}
L_n^{c\pm p -1}(zx)L_{m-s}^{c+s-1}(w\, x)dx\notag\\
&=\dfrac{(-1)^s(c)_m(c\pm p)_n\Gamma(c+j)}{m!\, n!\,\lambda^{c+j+s}}
\sum_{k=0}^n\dfrac{(-n)_k(c+j)_k}{(c\pm p)_k\,k!}\dfrac{(c+j+k)_s(-m)_s}{(c)_s}\left(\dfrac{z}{\lambda}\right)^k{}_2F_1\left(c+j+k+s,s-m;c+s;\dfrac{w}{\lambda}\right),\notag\\
&(m\geq s; c+j+s>0;~\lambda>0;~c+s,~c\pm p\neq 0,-1,\dots;p=0,1,\dots).
\end{align}
and further differentiation of both sides $\mu$-times, with respect to $z$, 
\begin{align}\label{integral81}
\int_0^\infty& x^{c+j+s+\mu-1}\,e^{-\lambda x}\,
L_{n-\mu}^{c\pm p+\mu -1}(zx)\,L_{m-s}^{c+s-1}(w\, x)\,dx\notag\\
&=\dfrac{(-m)_s (c)_m(c\pm p)_n\Gamma(c+j)}{(-1)^s\,m!\, n!\,z^\mu\,\lambda^{c+j+s}(c)_s}
\sum_{k=\mu}^n\dfrac{(-k)_\mu(c+j+k)_s(-n)_k(c+j)_k}{(c\pm p)_k\,k!}\left(\dfrac{z}{\lambda}\right)^k{}_2F_1(c+j+k+s,s-m;c+s;\dfrac{w}{\lambda}),\notag\\
&(s\leq m; ~\mu\leq n; ~c+j+s+\mu>0;~\lambda>0;~c,~c\pm p\neq 0,-1,\dots;p=0,1,\dots).
\end{align}
If $w=\lambda$, equation \eqref{integral79} reads
\begin{align}\label{integral82}
\int_0^\infty x^{c+j-1}e^{-\lambda x}
L_n^{c\pm p -1}(zx)L_m^{c-1}(\lambda\, x)dx
&=\dfrac{(-j)_m(c\pm p)_n\Gamma(c+j)}{m!\, n!\,\lambda^{c+j}}
{}_3F_2\left(-n,c+j,1+j;c\pm p,1+j-m;\dfrac{z}{\lambda}\right),\notag\\
&(j\geq m;~c+j>0;~\lambda>0;~c,~c\pm p\neq 0,-1,\dots;p=0,1,\dots).
\end{align}
and if $p=0$ it yields
\begin{align}\label{integral83}
\int_0^\infty x^{c+j-1}e^{-\lambda x}
L_n^{c+j -1}(zx)L_m^{c-1}(\lambda\, x)dx
&=\dfrac{(-j)_m\Gamma(c+j+n)}{m!\, n!\,\lambda^{c+j}}
{}_2F_1\left(-n,1+j;1+j-m;\dfrac{z}{\lambda}\right),\notag\\
&(j\geq m;~c+j>0;~\lambda>0).
\end{align}
and by taken limit of both sides as $j\rightarrow 0$
\begin{align}\label{integral84}
\int_0^\infty x^{c-1}e^{-\lambda x}
L_n^{c -1}(zx)L_m^{c-1}(\lambda\, x)dx
&=\dfrac{(-1)^{m}\,z^m\Gamma(c+n)(-n)_{m}}{m!\,n!\,\lambda^{c+n}\,(\lambda-z)^{m-n}}, \quad (n\geq m;~c>0;\lambda>0;|z|<\lambda),
\end{align}
thus
\begin{align}\label{integral85}
\int_0^\infty x^{c-1}e^{-\lambda x}
L_n^{c-1}(\lambda x)L_m^{c-1}(\lambda\, x)dx
&=\dfrac{(c)_n\Gamma(c)}{m!\, \lambda^{c}}
\delta_{m,n},\qquad
(c>0;~\lambda>0;n,m=0,1,\dots).
\end{align}
By means of 
$H_{2n}(\sqrt{z})=(-1)^n\,(2n)!{}_1F_1(-n;0.5;z)/n!$, it follows using \eqref{integral74} that
\begin{align}\label{integral86}
\int_0^\infty x^{j-\frac12}e^{-\lambda x}&\, L_{n}^{\pm p-\frac12}(zx)\, H_{2n}(\sqrt{wx})\, dx
=\dfrac{(-1)^n(2n)!(\pm p+\frac12)_n\Gamma(j+\frac12)}{(n!)^2\lambda^{j+\frac12}}\notag\\
&\times \sum_{k=0}^{n}\dfrac{(-n)_k(j+\frac12)_k}{(\pm p+\frac12)_k\,k!}\left(\dfrac{z}{\lambda}\right)^k{}_2F_1\left(j+k+\dfrac12,-n;\dfrac12;\dfrac{w}{\lambda}\right),\quad (j> 1/2;p,n=0,1,\dots),
\end{align}
from which it follows
\begin{align}\label{integral87}
\int_0^\infty x^{j-\frac12}e^{-\lambda x}&\, L_{n}^{\pm p-\frac12}(zx)\, H_{2n}(\sqrt{\lambda x})\, dx
=\dfrac{(-1)^n(2n)!(-j)_n(\pm p+\frac12)_n\Gamma(j+\frac12)}{(n!)^2\,(\frac12)_n\,\lambda^{j+\frac12}}\notag\\
&\times {}_3F_2\left(j+\frac12,1+j,-n;1+j-n,\pm p+\frac12;\dfrac{z}{\lambda}\right),\quad (j>1/2;\lambda>0;~~p,n=0,1,\dots)
\end{align}
However, by means of 
\begin{align}\label{def88}
\lim_{j\rightarrow 0}\,(-j)_n\,\,{}_3F_2\left(-n,\frac12+j,j+1;\pm p+\frac12,j+1-n;1\right)=\dfrac{(2n)!}{4^n(\pm p+\frac12)_n}
\end{align}
it easily follows that
\begin{align}\label{integral89}
\int_0^\infty x^{-\frac12}e^{-\lambda x}\, L_{n}^{\pm p-\frac12}(\lambda x)\, H_{2n}(\sqrt{\lambda x})\, dx
&=\dfrac{(-1)^n((2n)!)^2\sqrt{\pi}}{4^n(n!)^2\,(\frac12)_n\,\sqrt{\lambda}},\quad (\lambda>0, p: ~arbitrary).
\end{align}
Note also, if $c=1/2$ and $p=1$, it easily follows
\begin{align}\label{integral90}
\int_0^\infty x^{j-1}&e^{-\lambda x} H_{2m}(\sqrt{wx}) H_{2n+1}(\sqrt{z\,x})
\,dx\notag\\
&=(-1)^{m+n}\dfrac{(2m)!}{m!}\dfrac{2\sqrt{z}(2n+1)!}{n!}\dfrac{\Gamma(\frac12+j)}{\lambda^{\frac12+j}}\sum_{k=0}^m\dfrac{(\frac12+j)_k(-m)_k}{(\frac12)_k\,k!}\left(\dfrac{w}{\lambda}\right)^k{}_2F_1(-n,\frac12+j+k;\frac32;\dfrac{z}{\lambda}),\notag\\
&(j>0;~\lambda>0;~m,n=0,1,\dots),
\end{align}
For $j>n$ and $z=\lambda$
\begin{align}\label{integral91}
\int_0^\infty x^{j-1}&e^{-\lambda x} H_{2m}(\sqrt{wx}) H_{2n+1}(\sqrt{\lambda\,x})
\,dx=(-1)^{m+n}\dfrac{(2m)!}{m!}\dfrac{2(2n+1)!}{n!}\dfrac{\Gamma(\frac12+j)}{\lambda^{j}}\dfrac{(1-j)_n}{(\frac32)_n}
{}_3F_2(j+\dfrac12,j,-m;\dfrac12,j-n;\dfrac{w}{\lambda})
\end{align}
and for $c=1/2$ and $p=0$
\begin{align}\label{integral92}
\int_0^\infty x^{j-\frac12}&e^{-\lambda x}
H_{2m}(\sqrt{wx}) H_{2n}(\sqrt{z\,x})
\,dx=(-1)^{m+n}\dfrac{(2m)!}{m!}\dfrac{(2n)!}{n!}\dfrac{\Gamma(\frac12+j)}{\lambda^{\frac12+j}}\sum_{k=0}^m\dfrac{(\frac12+j)_k(-m)_k}{(\frac12)_k\,k!}\left(\dfrac{w}{\lambda}\right)^k{}_2F_1(-n,\frac12+j+k;\frac12;\dfrac{z}{\lambda})\end{align}
We may remark that all the above results involving ${}_1F_1$ can be rewritten in the representation using the
Whittaker function because of the following relationship:

$$
{}_1F_1(a,b,z)=e^{z/2}z^{-b/2}M_{(b-2a)/2,(b-1)/2}(z)
$$

\section*{Acknowledgements}

\noindent  This work was 
supported by the grant No. GP249577 from the Natural 
Sciences and Engineering Research Council of Canada.

\section{Appendix}
\noindent In this appendix we summarize some recurrence relations of the Gordon's integral that follow using the contiguous relations of the confluent hypergeometric functions:
\begin{align}
\operatorname {J}_{c+1}^{j(\pm p)}(b+1,b';\lambda,w,z)&=\dfrac{c}{w}\left[\operatorname {J}_{c}^{j(\pm p)}(b+1,b';\lambda,w,z)-\operatorname {J}_{c}^{j(\pm p)}(b,b';\lambda,w,z)\right].\label{A1}\\ \notag\\
\operatorname {J}_{c}^{j(\pm p)}(b+1,b';\lambda,w,z)&=\dfrac{c-b}{b}\operatorname {J}_{c}^{j(\pm p)}(b-1,b';\lambda,w,z)+\dfrac{w}{b}\operatorname {J}_{c}^{(j+1)(\pm p)}(b,b';\lambda,w,z)\notag\\
&+\dfrac{2b-c}{b}\operatorname {J}_{c}^{j(\pm p)}(b,b';\lambda,w,z),\label{A2}\\ \notag\\
\operatorname {J}_{c+1}^{j(\pm p)}(b+1,b';\lambda,w,z)&=\dfrac{(c+1-k)_k}{w^k}\sum_{m=0}^k\dfrac{(-1)^m\,k!}{m!\,(k-m)!}\operatorname {J}_{c+1-k}^{j(\pm p+k)}(b+1-m,b';\lambda,w,z)\label{A3}\\\notag\\
\operatorname {J}_{c}^{j(\pm p)}(b,b';\lambda,w,z)&=\dfrac{b}{c}\operatorname {J}_{c+1}^{(j-1)(\pm p-1)}(b+1,b';\lambda,w,z)-\dfrac{b-c}{c}\operatorname {J}_{c+1}^{(j-1)(\pm p-1)}(b,b';\lambda,w,z),\label{A4}\\\notag\\
\operatorname {J}_{c}^{j(\pm p)}(b+1,b';\lambda,w,z)&=\dfrac{b}{c}\operatorname {J}_{c+1}^{(j+1)(\pm p-1)}(b+1,b';\lambda,w,z)+\dfrac{w}{c}\operatorname {J}_{c+1}^{j(\pm p-1)}(b+1,b';\lambda,w,z)\notag\\
&-\dfrac{b-c}{c}\operatorname {J}_{c+1}^{(j+1)(\pm p-1)}(b,b';\lambda,w,z),\label{A5}\\ \notag\\
\operatorname {J}_{c}^{j(\pm p)}(b+1,b';\lambda,w,z)&=
\operatorname {J}_{c}^{j(\pm p)}(b,b';\lambda,w,z)
+\dfrac{w}{b}\operatorname {J}_{c}^{(j+1)(\pm p)}(b,b';\lambda,w,z)
-\dfrac{w(b-c)}{cb}\operatorname {J}_{c+1}^{j(\pm p-1)}(b,b';\lambda,w,z),\label{A6}\\ \notag\\
\operatorname {J}_{c-1}^{j(\pm p)}(b,b';\lambda,w,z)&=
\operatorname {J}_{c}^{(j-1)(\pm p)}(b,b';\lambda,w,z)
+\dfrac{w}{c-1}\operatorname {J}_{c}^{j(\pm p-1)}(b,b';\lambda,w,z)+\dfrac{w(b-c)}{c(1-c)}\operatorname {J}_{c+1}^{(j-1)(\pm p-2)}(b,b';\lambda,w,z),\label{A7}\\\notag\\
\operatorname {J}_{c}^{j0}(b,b';\lambda,w,z)&=\dfrac{b1-c}{z}
\operatorname {J}_{c}^{(j-1)0}(b,b'+1;\lambda,w,z)
+\dfrac{b'-c}{z}\operatorname {J}_{c}^{(j-1)0}(b,b'-1;\lambda,w,z)
+\dfrac{c-2b'}{z}\operatorname {J}_{c}^{(j-1)0}(b,b';\lambda,w,z).\label{A8}
\end{align}

\end{document}